\documentclass[11pt]{article}
\usepackage[english]{babel}
\usepackage[all]{xy}
\usepackage{fancyhdr}
\usepackage{setspace, amsmath, amsthm, amssymb, amsfonts, amscd, epic, graphicx, ulem, dsfont}
\usepackage{multirow}
\usepackage{amsthm}
\usepackage[T1]{fontenc}
\usepackage{answers}

\makeatletter

\@addtoreset{equation}{section}
\makeatother

\title{On torse-forming vector fields and biharmonic hypersurfaces in Riemannian manifolds}

\author{Ahmed Mohammed Cherif\footnote{University Mustapha Stambouli Mascara, Faculty of Exact Sciences, Mascara 29000, Algeria. Email: a.mohammedcherif@univ-mascara.dz}}
\date{}
\begin{document}
\maketitle

\begin{abstract}
In this paper, we give some properties of biharmonic hypersurface in Riemannian manifold has a torse-forming vector field.
\end{abstract}

\begin{flushleft}
Subjclass: 53C20; 58E20.\\
\end{flushleft}

\maketitle

\section{Introduction}
The energy functional of a smooth map $\varphi:(M,g)\longrightarrow(N,h)$ between two Riemannian manifolds is defined by
\begin{equation}\label{eq1.1}
    E(\varphi)=\frac{1}{2}\int_D|d\varphi|^2 v^g,
\end{equation}
where $D$ is compact domain of $M$,
$|d\varphi|$ is the Hilbert-Schmidt norm of the differential $d\varphi$, and $v^g$ is the volume element on $(M,g)$. A map $\varphi$
is called harmonic  if it is a critical point of the energy functional (\ref{eq1.1}).  The Euler Lagrange equation
associated to (\ref{eq1.1}) is given by (see \cite{BW,ES,YX})
\begin{equation}\label{eq1.2}
\tau(\varphi)=\operatorname{trace}\nabla d\varphi=\sum_{i=1}^m\nabla^\varphi _{e_i} d\varphi(e_i)-\sum_{i=1}^md\varphi(\nabla^M _{e_i}e_i)=0,
\end{equation}
where $\{e_i\}_{i=1}^m$ is a local orthonormal frame field on $(M,g)$, $\nabla^{M}$ is the Levi-Civita connection of $(M,g)$,
$\nabla^{\varphi}$ denote the pull-back connection on $\varphi^{-1}TN$, and $m$ is the dimension of $M$.
A natural generalization of harmonic maps is given by integrating the square of the norm of the tension field. More
precisely, the bienergy functional of a map $\varphi\in C^\infty(M,N)$ is defined by
\begin{equation}\label{eq1.3}
    E_2(\varphi)=\frac{1}{2}\int_D|\tau(\varphi)|^2 v^g.
\end{equation}
A map $\varphi\in C^\infty(M,N)$ is called biharmonic if it is a critical point of the bienergy functional, that is, if it is a solution of
the Euler Lagrange equation associated to (\ref{eq1.3})
\begin{eqnarray}\label{eq1.4}
% \nonumber to remove numbering (before each equation)
\tau_2(\varphi)
   &=&\nonumber -\operatorname{trace}R^N(\tau(\varphi),d\varphi)d\varphi-\operatorname{trace }\,(\nabla^{\varphi})^2 \tau(\varphi)  \\
   &=&\nonumber  -\sum_{i=1}^mR^N(\tau(\varphi),d\varphi(e_i))d\varphi(e_i)-\sum_{i=1}^m\nabla^{\varphi}_{e_i}\nabla^{\varphi}_{e_i}\tau(\varphi)\\
   &&+\sum_{i=1}^m\nabla^{\varphi}_{\nabla^M_{e_i}e_i}\tau(\varphi)=0,
\end{eqnarray}
where $R^N$ is the curvature tensor of $(N,h)$ defined by
\begin{equation*}
    R^N(X,Y)Z=\nabla^N_X \nabla^N_Y Z-\nabla^N_Y \nabla^N_X Z-\nabla^N_{[X,Y]}Z,
\end{equation*}
where $\nabla^N$ is the Levi-Civita connection of $(N,h)$ and $X,Y,Z\in\Gamma(TN)$
(see \cite{Jiang,YX}). Clearly, it follows from (\ref{eq1.4}) that any harmonic map is biharmonic and we call
those non-harmonic biharmonic maps proper biharmonic maps.\\
Let $M$ be a submanifold in $(N,\langle,\rangle)$ of dimension $m$, $ \mathbf{i}: M\hookrightarrow (N,\langle,\rangle)$  the canonical
inclusion, and let $\{e_i\}_{i=1}^m$ be a local orthonormal frame field with respect to induced Riemannian metric $g$ on $M$ by $\langle,\rangle$.
We denote by $\overline{\nabla}$ (resp. $\nabla$) the Levi-Civita connection of $(N,\langle,\rangle)$ (resp. of $(M,g)$),
by $\overline{\operatorname{grad}}$ (resp. $\operatorname{grad}$) the gradient operator on $(N,\langle,\rangle)$ (resp. on $(M,g)$),
by $B$ the second fundamental form of the submanifold $(M,g)$,
and by $H$ the mean curvature vector field of $(M,g)$
(see \cite{BW,ON}). The submanifold $(M,g)$ is called a harmonic (resp. biharmonic) submanifold in $(N,h)$
if $\tau(\mathbf{i})=0$ (resp. $\tau_2(\mathbf{i})=0)$. The expressions assumed by the tension and bitension fields are given by
\begin{equation}\label{eq1.5}
    \tau(\mathbf{i})=mH,\quad
    \tau_2(\mathbf{i})=-m\sum_{i=1}^m\left\{\overline{R}(H,e_i)e_i+
    \overline{\nabla}_{e_i}\overline{\nabla}_{e_i}H-\overline{\nabla}_{\nabla_{e_i}e_i}H\right\},
\end{equation}
where $\overline{R}$ is the curvature tensor of $(N,h)$. In \cite{OU}, Ye-Lin Ou proved that
a hypersurface $(M,g)$ in a Riemannian manifold $(N,\langle ,\rangle)$ with mean curvature vector field $H= f \eta $, that is the dimension of $N$ is $m+1$, is biharmonic if and only if
\begin{equation}\label{S}
\left\{
\begin{array}{lll}
-\Delta(f)  + f |A|^2 -f \operatorname{\overline{Ric}}(\eta , \eta )  &=& 0; \\\\
 2A(\operatorname{grad} f)+ mf \operatorname{grad}f -2 f (\operatorname{\overline{Ricci}} \eta)^\top  &=& 0,
\end{array}
\right.
\end{equation}
where $\operatorname{\overline{Ric}} $ (resp. $\operatorname{\overline{Ricci}}$) is the Ricci curvature (resp. Ricci tensor) of $(N,\langle ,\rangle )$,
$f$ denote the mean curvature function of $(M,g)$, and $A$ the shape operator with respect to the unit normal vector field $\eta$.\\
Let $(N,\langle ,\rangle)$ be a Riemannian manifold admits a torse-forming vector field $P$, that is $P$ satisfies the following formula
\begin{equation}\label{tf}
\overline{\nabla}_X P=\mu \,X+\omega(X)P,\quad\forall X\in\Gamma(TN),
\end{equation}
for some smooth function $\mu$ and $1$-form on $N$. The $1$-form $\omega$ is called the generating form
and the function $\mu$ is called the conformal scalar.
Let $(M,g)$ be a hypersurface in $(N,\langle ,\rangle)$.
We consider the following decomposition of the torse-forming vector field
$$P=\phi\,\eta+V,$$ where $V$ denote the tangential component of $P$ and $\phi=\langle P,\eta\rangle$.


\begin{thebibliography}{0}
\bibitem{BW} P. Baird and J. C. Wood,  {\it Harmonic morphisms between Riemannain manifolds}, Clarendon Press Oxford (2003).

\bibitem{CM} R. Caddeo and S. Montaldo, {\it C. Oniciuc. Biharmonic submanifolds in spheres}, Israel J. Math., {\bf130} (2002), 109-123.


\bibitem{chen2} B-Y. Chen, {\it Some open problems and conjectures on submanifolds of finite type}, Soochow J. Math., {\bf17} (1991), 169-188.


\bibitem{chen3} B.Y. Chen, {\it A report on submanifolds of finite type}, Soochow J. Math., {\bf 22} (1996), 117-337.


\bibitem{ES} J. Eells and J. H. Sampson, {\it Harmonic mappings of Riemannian manifolds}, { Amer. J. Math.,} {\bf 86} (1964), 109-160.

\bibitem{Jiang} G. Y. Jiang, {\it $2$-Harmonic maps between Riemannian manifolds}, Annals of Math., China, {\bf 7A(4)} (1986), 389-402.


\bibitem{kenmotsu} K. Kenmotsu, {\it A class of almost contact Riemannian manifolds}, Tohoku Math.J., {\bf 24} (1972), 93-103.


\bibitem{cd} A. Mohammed Cherif and M. Djaa, {\it Harmonic maps and torse-forming vector
fields}, International electronic journal of geometry, {\bf 13}(1) (2020), 87-93.


\bibitem{NU} N. Nakauchi and H. Urakawa, {\it Biharmonic hypersurfaces in a Riemannian
manifold with non-positive Ricci curvature}, Ann. Global Anal. Geom., {\bf40} (2011), 125-131.

\bibitem{ON} O'Neil, {\it Semi-Riemannian Geometry}, Academic Press, New York (1983).


\bibitem{OU} Y-L. Ou, {\it Biharmonic hypersurfaces in Riemannian manifolds}, Pacific Journal of Mathematics, {\bf 248}(1) (2010), 217-232.

\bibitem{YX} Y. Xin, {\it Geometry of harmonic maps}, Fudan University (1996).



\end{thebibliography}
\end{document}